\documentclass[12pt]{article}%
\usepackage{eurosym}
\usepackage{latexsym}
\usepackage{graphicx}
\usepackage{enumerate}
\usepackage{amsthm}
\usepackage{cite}
\usepackage{amsmath}
\usepackage{amsfonts}
\usepackage{amssymb}
\usepackage{graphicx}
\usepackage{psfrag}
\usepackage{subfig}
\usepackage{color}
\usepackage[normalem]{ulem}%
\providecommand{\U}[1]{\protect\rule{.1in}{.1in}}
\graphicspath{{c:/swp55/Graphics/}}
\newtheorem{theorem}{Theorem}[section]

\newtheorem{lemma}[theorem]{Lemma}
\newtheorem{notation}[theorem]{Notation}

\newtheorem{remark}[theorem]{Remark}

\newcommand{\BIGOP}[1]{\mathop{\mathchoice{\raise-0.22em\hbox{\huge
$#1$}} {\raise-0.05em\hbox{\Large $#1$}}{\hbox{\large $#1$}}{#1}}}
\makeatletter\@addtoreset{equation}{section}\makeatother
\evensidemargin=\oddsidemargin
\newdimen\dummy
\dummy=\oddsidemargin
\makeatletter
\g@addto@macro\@floatboxreset\centering
\makeatother
\begin{document}

\title{A Posteriori Modelling-Discretization Error Estimate for Elliptic Problems
with $L^{\infty}$-Coefficients}
\author{M. Weymuth\thanks{Institut f\"{u}r Mathematik, Universit\"{a}t Z\"{u}rich,
Winterthurerstrasse 190, CH-8057 Z\"{u}rich, Switzerland, e-mail:
\texttt{monika.weymuth@math.uzh.ch}}
\and S. Sauter\thanks{Institut f\"{u}r Mathematik, Universit\"{a}t Z\"{u}rich,
Winterthurerstrasse 190, CH-8057 Z\"{u}rich, Switzerland, e-mail:
\texttt{stas@math.uzh.ch}}
\and {S. Repin}\thanks{V.A. Steklov Institute of Mathematics, Fontanka 27, 191 011
St. Petersburg, Russia, e-mail: \quad\texttt{repin@pdmi.ras.ru}; } }
\maketitle

\begin{abstract}
We consider elliptic problems with complicated, discontinuous diffusion tensor
$A_{\scriptscriptstyle 0} $. One of the standard approaches to numerically
treat such problems is to simplify the coefficient by some approximation, say
$A_{\varepsilon}$, and to use standard finite elements. In \cite{Repin2012} a
combined modelling-discretization strategy has been proposed which estimates
the discretization and modelling errors by a posteriori estimates of
functional type. This strategy allows to balance these two errors in a problem
adapted way. However, the estimate of the modelling error is derived under the
assumption that the difference $A_{\scriptscriptstyle 0} -A_{\varepsilon}$ is
bounded in the $L^{\infty}$-norm, which requires that the approximation of the
coefficient matches the discontinuities of the original coefficient. Therefore
this theory is not appropriate for applications with discontinuous
coefficients along \textit{complicated, curved} interfaces. Based on bounds
for $A_{\scriptscriptstyle 0} -A_{\varepsilon}$ in an $L^{q}$-norm with
$q<\infty$ we generalize the combined modelling-discretization strategy to a
larger class of coefficients.

\end{abstract}

\section{Introduction}

We consider elliptic boundary value problems with complicated, discontinuous
diffusion tensor. As a model problem we choose the diffusion equation
$-\operatorname{div}(A_{\scriptscriptstyle 0} \nabla u)=f$ in a two- or
three-dimensional bounded domain $\Omega$ with homogeneous Dirichlet boundary
conditions and $A_{\scriptscriptstyle 0} \in L^{\infty}(\Omega,\mathbb{R}%
_{\operatorname*{sym}}^{d\times d})$ is symmetric and positive definite. Our
emphasis is on diffusion matrices $A_{\scriptscriptstyle 0} $ containing a
large number of different scales which we allow to be highly non-uniformly
distributed over the domain.

It is well-known that for such problems standard single scale numerical
methods such as standard finite element methods are not efficient, since one
needs to solve the problem on a sufficiently fine mesh which resolves all the
fine-scale behavior of the coefficient. This is usually too costly, especially
for three-dimensional problems. Essentially there are two different approaches
to overcome this difficulty: One is to design (non-polynomial) generalized
finite element methods where the characteristic behavior of the solution is
reflected by the shape of the basis functions. This approach has been
investigated by many researchers (see, e.g., \cite{Babuska1994},
\cite{Babuska1997}, \cite{Babuska2004}, \cite{Melenk1996},
\cite{Strouboulis2000_1,Strouboulis2000_2,Strouboulis2001}). In this paper we
follow the second approach which tries to simplify the diffusion coefficient
by some approximation and then employs standard finite elements. Standard
methods for simplifying the coefficients are based, e.g., on homogenization
methods for periodic structures (see. e.g., \cite{jikov94}, \cite{Cioranescu},
\cite{Bensoussan}), or on different upscaling techniques e.g. \cite{E_Enquist}%
, \cite{PreRumSau2011}, \cite{Repin2012}.

In many applications a numerical solution with only moderate guaranteed
accuracy is required. For such problems, the \textit{combined
modelling-discretization strategy} has been proposed in \cite{Repin2012}. This
approach consists of two basic steps. In a first step the diffusion
coefficient $A_{\scriptscriptstyle 0}$ is replaced by a simpler coefficient
$A_{\varepsilon}$ and the simplified model is discretized and solved on a
rather coarse mesh. In the second step the discretization and modelling errors
are controlled using some a posteriori estimates. The total error is bounded
by the sum of the discretization and modelling errors which are both explicit
and computable. If the total error is larger than a given tolerance, then
either the mesh should be refined (if the discretization error dominates) or
the coefficient has to be modelled more accurately (if the modelling error
dominates). Thus the discretization and modelling errors can be balanced in a
problem-adapted way.

The error estimates in \cite{Repin2012} are derived by purely functional
methods without requiring specific information on the approximating subspace
and the numerical method used. Consequently the estimates contain no mesh
dependent constants and are valid for any conforming approximation from the
respective energy space.

However, the modelling errors arising due to the simplification of the
coefficients contain the term $\left\vert \left\vert \left\vert
A_{\scriptscriptstyle 0} -A_{\varepsilon}\right\vert \right\vert \right\vert
_{\infty,\Omega}$ (cf. \cite{Repin2012}) and one has to assume that the
approximation $A_{\varepsilon}$ matches the discontinuities of
$A_{\scriptscriptstyle 0}$ in order to ensure that the term $\left\vert
\left\vert \left\vert A_{\scriptscriptstyle 0} -A_{\varepsilon}\right\vert
\right\vert \right\vert _{\infty,\Omega}$ becomes small. Since in many
applications the discontinuities jump on curved or cracked interfaces, they
cannot be captured exactly by the finite element mesh. Hence, this smallness
assumption is not suitable for the analysis of numerical methods for problems
with discontinuous coefficients along complicated, curved interfaces.

This problem has been addressed in \cite{Nochetto2013} in the context of
adaptive finite element methods. Based on a perturbation theory the term
$\left\vert \left\vert \left\vert A_{\scriptscriptstyle 0} -A_{\varepsilon
}\right\vert \right\vert \right\vert _{\infty,\Omega}$ is replaced by
$\left\vert \left\vert \left\vert A_{\scriptscriptstyle 0} -A_{\varepsilon
}\right\vert \right\vert \right\vert _{q,\Omega}$ with $q:=2p/(p-2)$ for some
$p\geq2$. The advantage of this approach is that $A_{\varepsilon}$ does not
have to match the discontinuities of $A_{\scriptscriptstyle 0} $ exactly.
However, one needs more regularity on the solution, namely $\nabla u\in
L^{p}(\Omega)$ for some $p>2$. This also requires additional assumptions on
the right-hand side $f$. These requirements are quite mild and are satisfied
in many applications.

The goal of our paper is to generalize the modelling-discretization strategy
developed in \cite{Repin2012} to a larger class of coefficients. Based on the
theory presented in \cite{Nochetto2013} which is based on results by
\cite{Meyers1963} we bound the modelling error by a term depending on
$\left\vert \left\vert \left\vert A_{\scriptscriptstyle 0} -A_{\varepsilon
}\right\vert \right\vert \right\vert _{q,\Omega}$ for some $q>2$. Consequently
the assumption on $A_{\varepsilon}$ can be weakened and the strategy can also
be applied to problems where $A_{\scriptscriptstyle 0} $ has discontinuities
which are unknown or lie along curves and surfaces. Moreover we have
guaranteed, computable upper bounds of the error.

Due to our new approach the bound of the modelling error also depends on a
regularity constant $C_{\operatorname*{reg},A_{\varepsilon}}$ with respect to
some Sobolev space $W^{k,p}\left(  \Omega\right)  $ as it appears in the
inequality
\[
\Vert\nabla u_{\varepsilon}\Vert_{p,\Omega}\leq C_{\operatorname*{reg}%
,A_{\varepsilon}}\Vert F\Vert_{-1,p,\Omega}%
\]
for some $p>2$, where $u_{\varepsilon}$ is the exact solution of the
simplified problem and $F$ is a linear functional generated by the right-hand
side of the equation. The constant $C_{\operatorname*{reg},A_{\varepsilon}}$
depends only on $p$, $A_{\varepsilon}$ and $\Omega$. In \cite{Meyers1963,
Nochetto2013}, it is shown (by perturbation arguments) that
$C_{\operatorname*{reg},A_{\varepsilon}}$ can be expressed in terms of the
constant $C_{\operatorname*{reg},I}$ which corresponds to the Laplace
operator. Therefore, we need to derive computable upper bounds for
$C_{p}:=C_{\operatorname*{reg},I}$.

The paper is structured as follows. In Section \ref{sec:Problem}, we first
formulate the model problem and the conditions on the coefficient. Then, in
Section \ref{sec:discretization} we present new error estimates for the
modelling-discretization strategy introduced in \cite{Repin2012}. These
estimates are based on the theory developed in \cite{Meyers1963, Nochetto2013}
and require that $\nabla u_{\varepsilon}\in L^{p}(\Omega)$ for some $p>2$.
Section \ref{sec:regularity} is devoted to present some $L^{p}$-bounds for the
gradient of the solution of diffusion problems with $L^{\infty}$-coefficient.
These bounds only depend on the size of the jumps in the coefficient. Finally
in Section \ref{sec:constant} we present an explicit computable estimate of
$C_{p}$ for the full space problem. This bound depends on $p$ and the
dimension $d$.

\section{Notation and Problem Statement}

\label{sec:Problem}

Throughout the paper it is assumed that $\Omega$ is a bounded domain in
$\mathbb{R}^{d}$ ($d=2,3$) with $C^{1}$ boundary. By $\langle\cdot
,\cdot\rangle$ we denote the usual Euclidean scalar product on $\mathbb{R}%
^{d}$. For $1\leq p\leq\infty$, $\Vert\cdot\Vert_{\ell^{p}}$ denotes the
discrete $\ell^{p}$-norm in $\mathbb{R}^{d}$. If $p=2$, then we use
$\left\Vert \cdot\right\Vert $ instead of $\left\Vert \cdot\right\Vert
_{\ell^{2}}$. By $\left(  \cdot,\cdot\right)  $ we denote the $L^{2}$-scalar
product. The Sobolev space of real-valued functions in $L^{2}\left(
\Omega\right)  $ with gradients in $L^{2}(\Omega)$ is denoted by $H^{1}%
(\Omega)$ (and $\|\cdot\|_{1,2,\Omega}$ is the respective norm). A subspace of
$H^{1}(\Omega)$ containing the functions vanishing on the boundary is denoted
by $H_{0}^{1}(\Omega)$. For any $p\in[1,\infty]$, the adjoint number
$p^{\prime}$ is defined by the relation $\frac{1}{p}+\frac{1}{p^{\prime}}=1$.
Analogously, for $p\in\left[  2,\infty\right]  $, the number $p^{\prime\prime
}\in\left[  1,\infty\right]  $ satisfies the relation $\frac{2}{p}+\frac
{1}{p^{\prime\prime}}=1$.

Throughout the paper $\Vert\cdot\Vert_{p,\Omega}$ denotes the norm of
$L^{p}(\Omega)$. We use the usual notation $W^{1,p}(\Omega)$ for the Sobolev
spaces of functions, which generalized derivatives belong to the space
$L^{p}(\Omega)$. This space is supplied with the standard norm $\|
\cdot\|_{1,p,\Omega}$. The space of functions denoted by $W_{0}^{1,p}(\Omega)$
is the closure of $C_{0}^{\infty}(\Omega)$ with respect to the norm $\|
\cdot\|_{1,p,\Omega}$. We also use the space $W^{-1,p}\left(  \Omega\right)
:=(W_{0}^{1,p^{\prime}}\left(  \Omega\right)  )^{\prime}$ endowed with the
standard dual norm $\| \cdot\|_{-1,p,\Omega}$. For vector and matrix valued
functions, we use the same notation for the Lebesgue and Sobolev spaces as
well as for the corresponding norms. To explicitly indicate the dimension we
write $L^{p}\left(  \Omega,\mathbb{R}^{d}\right)  $ and $L^{p}\left(
\Omega,\mathbb{R}^{d\times d}\right)  $. For functions in $L^{2}\left(
\Omega,\mathbb{R}^{d}\right)  $ we set
\[
\Vert\cdot\Vert_{p,\Omega}:=\Vert\,\,\Vert\cdot\Vert_{\ell^{p}}\,\Vert
_{p,\Omega}.
\]
For $M\in L^{\infty}\left(  \Omega,\mathbb{R}_{\operatorname{sym}}^{d\times
d}\right)  $ and $p\geq2$, we introduce the function $m\in L^{\infty}\left(
\Omega\right)  $ by%
\[
m:=\sup_{\zeta\in\mathbb{R}^{d}\backslash\left\{  0\right\}  }\frac{\left\Vert
M\left(  \cdot\right)  \zeta\right\Vert _{\ell^{p^{\prime}}}}{\left\Vert
\zeta\right\Vert _{\ell^{p}}}\quad\text{and\quad} \left\vert \left\vert
\left\vert M\right\vert \right\vert \right\vert _{p^{\prime\prime},\Omega
}=\left\Vert m\right\Vert _{p^{\prime\prime},\Omega}.
\]
Notice that for $p=2$ we have $p^{\prime}=2$ and $p^{\prime\prime}=\infty$ so
that
\begin{equation}
\left\vert \left\vert \left\vert M\right\vert \right\vert \right\vert
_{\infty,\Omega}=\underset{x\in\Omega}{\operatorname*{ess}\sup}\left(
\sup_{\zeta\in\mathbb{R}^{d}\backslash\left\{  0\right\}  }\frac{\left\Vert
M\left(  x\right)  \zeta\right\Vert }{\left\Vert \zeta\right\Vert }\right)  .
\label{spectrum}%
\end{equation}
Also we use the space
\[
H(\Omega,\operatorname{div}):=\{y\in L^{2}(\Omega,\mathbb{R}^{d}%
)\mid\operatorname{div}y\in L^{2}(\Omega)\},
\]
which is a Hilbert space endowed with the scalar product
\[
( y,z)_{\operatorname{div}}:=\left(  y,z\right)  +\left(  \operatorname*{div}%
y,\operatorname*{div}z\right)
\]
and the norm $\Vert y\Vert_{\operatorname{div}}:= ( y,y)_{\operatorname{div}%
}^{1/2}$. For the functions in $L^{2}(\Omega,\mathbb{R}^{d})$, we also
introduce two equivalent norms (associated with the energy and complementary
energy)
\[
\Vert y\Vert_{A_{\scriptscriptstyle 0} }^{2}:=\left(  A_{\scriptscriptstyle 0}
y,y\right)  =\int\limits_{\Omega}\langle A_{\scriptscriptstyle 0}
y,y\rangle\qquad\text{and}\quad\Vert y\Vert_{A^{-1}_{0}}^{2}:=\left(
A^{-1}_{0} y,y\right)  ,
\]
where the matrix $A_{\scriptscriptstyle 0} \in L^{\infty}\left(
\Omega,\mathbb{R}_{\operatorname*{sym}}^{d\times d}\right)  $ is assumed to be
uniformly positive definite, i.e.%
\begin{equation}
0<\alpha_{\scriptscriptstyle0}:=\left\vert \left\vert \left\vert
A^{-1}_{\scriptscriptstyle 0}\right\vert \right\vert \right\vert
_{\infty,\Omega}^{-1}\leq\left\vert \left\vert \left\vert
A_{\scriptscriptstyle 0} \right\vert \right\vert \right\vert _{\infty,\Omega
}=:\beta_{\scriptscriptstyle0}<\infty. \label{defalphabeta}%
\end{equation}

Let $f$ be a given function in $L^{2}\left(  \Omega\right)  $. Consider the
following boundary value problem: Find $u\in H_{0}^{1}\left(  \Omega\right)  $
such that%
\begin{equation}
\int_{\Omega}\left\langle A_{\scriptscriptstyle 0} \nabla u,\nabla
v\right\rangle =\int_{\Omega}fv=:F(v)\qquad\forall v\in H_{0}^{1}\left(
\Omega\right)  . \label{mainproblem}%
\end{equation}

In view of (\ref{defalphabeta}), existence and uniqueness of the solution $u$
follows from the Lax-Milgram lemma.

We consider problems with complicated matrix $A_{\scriptscriptstyle0}$ (which
coefficients are complicated functions of $x$). In this case, direct
approximation of $u$ based upon standard numerical approaches may lead to high
computational costs. One way to obtain a reasonable approximation of $u$ with
minimal expenditures is to consider a simplified problem with a simpler matrix
$A_{\varepsilon}$. If the difference between $u$ and the respective solution
$u_{\varepsilon}$ is explicitly estimated and it is smaller than the desired
tolerance, then we can use the simplified problem instead of the problem
(\ref{mainproblem}). This idea leads to a set of simplified problems generated
by uniformly positive definite matrices $A_{\varepsilon}\in L^{\infty}%
(\Omega,\mathbb{R}_{\operatorname*{sym}}^{d\times d})$, where $\varepsilon$ is
a sequence of positive decreasing numbers (which is either finite, or infinite
tending to zero). Henceforth, we assume that there exist positive constants
$\underline{\alpha}$ and $\overline{\beta}$ such that for any $\varepsilon$
\begin{equation}
0<\underline{\alpha}\leq\alpha_{\varepsilon}:=\left\vert \left\vert \left\vert
A_{\varepsilon}^{-1}\right\vert \right\vert \right\vert _{\infty,\Omega}%
^{-1}\leq\left\vert \left\vert \left\vert A_{\varepsilon}\right\vert
\right\vert \right\vert _{\infty,\Omega}=:\beta_{\varepsilon}\leq
\overline{\beta}<\infty\label{simpl_coeff}%
\end{equation}
and that the collection of simplified matrices $A_{\varepsilon}$ satisfies the
condition
\begin{equation}
\Vert A_{\scriptscriptstyle0}-A_{\varepsilon}\Vert_{q,\Omega}\leq\varepsilon.
\label{A0Aepsilon}%
\end{equation}


\section{Discretization and Combined Error Majorant}

\label{sec:discretization} The function $u_{\varepsilon}\in H_{0}^{1}\left(
\Omega\right)  $ (generalized solution of the simplified problem) is defined
by the integral identity
\begin{equation}
\int_{\Omega}\left\langle A_{\varepsilon}\nabla u_{\varepsilon},\nabla
v\right\rangle =F(v)\qquad\forall v\in H_{0}^{1}\left(  \Omega\right)
.\label{simpl_problem}%
\end{equation}
The difference between $u$ and $u_{\varepsilon}$ is the \textit{modelling
error}
\[
E_{\operatorname*{mod}}^{\varepsilon}:=\Vert\nabla(u-u_{\varepsilon}%
)\Vert_{A_{\scriptscriptstyle0}}.
\]
In general $u_{\varepsilon}$ is unknown and instead we use a conforming
approximation $u_{\varepsilon,h}\in H_{0}^{1}\left(  \Omega\right)  $ computed
by some numerical method. In view of this, we must also consider the
\textit{discretization error}
\[
E_{\operatorname*{disc}}^{\varepsilon,h}:=\Vert\nabla(u_{\varepsilon
}-u_{\varepsilon,h})\Vert_{A_{\scriptscriptstyle0}}.
\]
At this point, we do not specify the method by which $u_{\varepsilon,h}$ is
found. In the framework of our approach it is not important because the
difference between $u_{\varepsilon}$ and $u_{\varepsilon,h}$ can be estimated
within the framework of a unified method that follows from a posteriori error
estimates of the functional type (see e.g. \cite{Repin2008} and the references
therein). In our case, the respective estimate has the form
\begin{align}
\Vert\nabla(u_{\varepsilon}-u_{\varepsilon,h})\Vert_{A_{\varepsilon}}^{2}\, &
\leq\,\mathcal{M}_{\Omega}^{2}(u_{\varepsilon,h},y,\gamma
)\nonumber\label{majorant}\\
&  :=(1+\gamma)\Vert A_{\varepsilon}\nabla u_{\varepsilon,h}-y\Vert
_{A_{\varepsilon}^{-1}}^{2}+\left(  1+\frac{1}{\gamma}\right)  C_{\Omega}%
^{2}\Vert\operatorname{div}y+f\Vert_{2,\Omega}^{2}.
\end{align}
The majorant $\mathcal{M}_{\Omega}^{2}(u_{\varepsilon,h},y,\gamma)$ contains a
vector-valued function $y\in H(\Omega,\operatorname{div})$ and an arbitrary
positive parameter $\gamma$. The constant depends on the geometry of $\Omega$,
namely, $C_{\Omega}:=C_{F\Omega}/\sqrt{\underline{\alpha}}$, where
\[
C_{F\Omega}:=\sup_{w\in H_{0}^{1}(\Omega)\backslash\{0\}}\frac{\Vert
w\Vert_{2,\Omega}}{\Vert\nabla w\Vert_{2,\Omega}}\,\leq\,\frac{\mathrm{diam}%
\,\Omega}{\sqrt{2}\pi}.
\]
In order to formulate the main result, we introduce the quantities
\[
\Upsilon(u_{\varepsilon,h},f,\theta(p,t)):=\left(  C_{\operatorname*{reg}%
,A_{\varepsilon}}\Vert f\Vert_{t,\Omega}+\Vert\nabla u_{\varepsilon,h}%
\Vert_{t,\Omega}\right)  ^{1-\theta(p,t)}{\alpha_{\varepsilon}^{-\frac{1}%
{2}\theta(p,t)}}%
\]
and
\[
\Theta(u_{\varepsilon,h}):=\left(  \mathcal{M}_{\Omega}(u_{\varepsilon
,h},y,\gamma)\right)  ^{\theta(p,t)}\Upsilon(u_{\varepsilon,h},f,\theta
(p,t))+\Vert\nabla u_{\varepsilon,h}\Vert_{p,\Omega}.
\]
Here $\theta(r,t):=\frac{2(t-r)}{r(t-2)}$ (for $2<r<t<\infty$) and
$C_{\operatorname*{reg},A_{\varepsilon}}$ is a constant in the inequality
\[
\Vert\nabla u_{\varepsilon}\Vert_{t,\Omega}\leq C_{\operatorname*{reg}%
,A_{\varepsilon}}\Vert F\Vert_{-1,t,\Omega}.
\]
In Section \ref{sec:regularity}, we show that $C_{\operatorname*{reg}%
,A_{\varepsilon}}$ depends only on the constant $C_{P}:=C_{\operatorname*{reg}%
,I}$ (associated with the Laplace operator) and on the amplitude of the jumps
in the coefficient $A_{\varepsilon}$ (cf. Theorem \ref{TheoReg}).

\begin{theorem}
\label{theo_error} Let $A_{\scriptscriptstyle0}\in L^{\infty}(\Omega
,\mathbb{R}_{\operatorname*{sym}}^{d\times d})$ satisfy \eqref{defalphabeta},
$f\in L^{P}\left(  \Omega\right)  $ for some $P\in(2,+\infty)$, and
$p\in(2,p^{\ast})$, where the function $p^{\ast}=p^{\ast}\left(
\frac{\underline{\alpha}}{\overline{\beta}},P\right)  $ is defined by
(\ref{defpast2}). Then
\[
\Vert\nabla(u-u_{\varepsilon,h})\Vert_{A_{\scriptscriptstyle0}}\leq
E_{\operatorname*{disc}}^{\varepsilon,h}+E_{\operatorname*{mod}}^{\varepsilon
},
\]
where
\begin{align}
&  E_{\operatorname*{disc}}^{\varepsilon,h}\leq\left\vert \left\vert
\left\vert D_{\varepsilon}\right\vert \right\vert \right\vert _{\infty,\Omega
}^{1/2}\mathcal{M}_{\Omega}(u_{\varepsilon,h},y,\gamma)\label{disc_error}\\
&  E_{\operatorname*{mod}}^{\varepsilon}\leq\left\vert \left\vert \left\vert
B_{\varepsilon}\right\vert \right\vert \right\vert _{p^{\prime\prime},\Omega
}^{1/2}\Theta(u_{\varepsilon,h}), \label{mod_error}%
\end{align}
$D_{\varepsilon}:=A_{\varepsilon}^{-1/2}A_{\scriptscriptstyle0}A_{\varepsilon
}^{-1/2}$, $B_{\varepsilon}:=(A_{\varepsilon}-A_{\scriptscriptstyle0}%
)A_{\scriptscriptstyle0}^{-1}(A_{\varepsilon}-A_{\scriptscriptstyle0})$, and
$p^{\prime\prime}$ is a number defined in Sect. \ref{sec:Problem}.
\end{theorem}

\proof
By the triangle inequality
\[
\Vert\nabla(u-u_{\varepsilon,h})\Vert_{A_{\scriptscriptstyle 0} }\leq
\Vert\nabla(u_{\varepsilon}-u_{\varepsilon,h})\Vert_{A_{\scriptscriptstyle 0}
}+\Vert\nabla(u-u_{\varepsilon})\Vert_{A_{\scriptscriptstyle 0} }%
=E_{\operatorname*{disc}}^{\varepsilon,h}+E_{\operatorname*{mod}}%
^{\varepsilon}.
\]
First we estimate the discretization error. It holds%
\begin{align*}
\left(  E_{\operatorname*{disc}}^{\varepsilon,h}\right)  ^{2}  &
=\int_{\Omega} \left\langle A_{\varepsilon}^{1/2}D_{\varepsilon}%
A_{\varepsilon}^{1/2}\nabla(u_{\varepsilon}-u_{\varepsilon,h}),\nabla
(u_{\varepsilon}-u_{\varepsilon,h})\right\rangle \\
&  =\int_{\Omega}\left\langle D_{\varepsilon}A_{\varepsilon}^{1/2}%
\nabla(u_{\varepsilon}-u_{\varepsilon,h}),A_{\varepsilon}^{1/2}\nabla
(u_{\varepsilon}-u_{\varepsilon,h})\right\rangle \\
&  \leq\left\vert \left\vert \left\vert D_{\varepsilon} \right\vert
\right\vert \right\vert _{\infty,\Omega}\int_{\Omega}\left\langle
A_{\varepsilon}\nabla(u_{\varepsilon}-u_{\varepsilon,h}),\nabla(u_{\varepsilon
}-u_{\varepsilon,h})\right\rangle \\
&  = \left\vert \left\vert \left\vert D_{\varepsilon} \right\vert \right\vert
\right\vert _{\infty,\Omega}\Vert\nabla(u_{\varepsilon}-u_{\varepsilon
,h})\Vert_{A_{\varepsilon}}^{2}.
\end{align*}
The last norm can be estimated by the error majorant and thus we obtain
\eqref{disc_error}. Next we will estimate the modelling part of the error.
Observe that
\[
0=\int\limits_{\Omega}\langle A_{\scriptscriptstyle 0} \nabla(u-u_{\varepsilon
}),\nabla v \rangle+\int\limits_{\Omega}\langle(A_{\scriptscriptstyle 0}
-A_{\varepsilon})\nabla u_{\varepsilon}, \nabla v\rangle\qquad\forall v\in
H_{0}^{1}(\Omega).
\]
We choose $v=u-u_{\varepsilon}$ and obtain
\begin{align*}
\left(  E_{\operatorname*{mod}}^{\varepsilon}\right)  ^{2}=\Vert
\nabla(u-u_{\varepsilon})\Vert_{A_{\scriptscriptstyle 0} }^{2}  &  =
\int\limits_{\Omega}\langle A_{\scriptscriptstyle 0} \nabla(u-u_{\varepsilon
}),\nabla(u-u_{\varepsilon})\rangle= \int\limits_{\Omega}\langle
(A_{\varepsilon}-A_{\scriptscriptstyle 0} )\nabla u_{\varepsilon},
\nabla(u-u_{\varepsilon}\rangle.
\end{align*}
Applying the Cauchy-Schwarz inequality yields%
\[
\Vert\nabla(u-u_{\varepsilon})\Vert_{A_{\scriptscriptstyle 0} }^{2}\leq\left(
\int_{\Omega}\langle B_{\varepsilon}\nabla u_{\varepsilon},\nabla
u_{\varepsilon}\rangle\right)  ^{1/2}\Vert\nabla(u-u_{\varepsilon}%
)\Vert_{A_{\scriptscriptstyle 0} }.
\]
Hence,
\[
\Vert\nabla(u-u_{\varepsilon})\Vert_{A_{\scriptscriptstyle 0} }^{2}\leq
\int_{\Omega}\langle B_{\varepsilon}\nabla u_{\varepsilon},\nabla
u_{\varepsilon}\rangle.
\]
We set $p^{\prime\prime}=p/(p-2)$, apply the triangle and generalized
H\"{o}lder inequalities, and obtain
\begin{align}
\Vert\nabla(u-u_{\varepsilon})\Vert_{A_{\scriptscriptstyle 0} }  &
\leq\left\vert \left\vert \left\vert B_{\varepsilon} \right\vert \right\vert
\right\vert _{p^{\prime\prime},\Omega}^{1/2}\Vert\nabla u_{\varepsilon}%
\Vert_{p,\Omega}\nonumber\label{eq:1}\\
&  \leq\left\vert \left\vert \left\vert B_{\varepsilon} \right\vert
\right\vert \right\vert _{p^{\prime\prime},\Omega}^{1/2} \left(  \Vert
\nabla(u_{\varepsilon}-u_{\varepsilon,h}) \Vert_{p,\Omega}+\Vert\nabla
u_{\varepsilon,h}\Vert_{p,\Omega}\right)  .
\end{align}
Now, we choose $t\in(p,p^{\ast}(\frac{\underline{\alpha}}{\overline{\beta}%
},P))$. By Lemma \ref{interpolation_estimate} and (\ref{majorant}), we find
that
\begin{align*}
\Vert\nabla(u_{\varepsilon}-u_{\varepsilon,h})\Vert_{p,\Omega} & \leq
\Vert\nabla(u_{\varepsilon}-u_{\varepsilon,h})\Vert_{2,\Omega}^{\theta
(p,t)}\Vert\nabla(u_{\varepsilon}-u_{\varepsilon,h})\Vert_{t,\Omega}%
^{1-\theta(p,t)}\\
& \leq\left(  \frac{1}{\sqrt{\alpha_{\varepsilon}}}\mathcal{M}_{\Omega
}(u_{\varepsilon,h},y,\gamma)\right)  ^{\theta(p,t)}\left(  \Vert\nabla
u_{\varepsilon}\Vert_{t,\Omega}+\Vert\nabla u_{\varepsilon,h}\Vert_{t,\Omega
}\right)  ^{1-\theta(p,t)}.
\end{align*}
Now, we apply Theorem \ref{TheoReg} (where the regularity constant
$C_{\operatorname*{reg},A_{\varepsilon}}$ is defined) and arrive at the
estimate
\begin{align}
\label{eq:2}\Vert\nabla(u_{\varepsilon}-u_{\varepsilon,h})\Vert_{p,\Omega}  &
\leq\left(  \frac{1}{\sqrt{\alpha_{\varepsilon}}}\mathcal{M}_{\Omega
}(u_{\varepsilon,h},y,\gamma)\right)  ^{\theta(p,t)}\left(
C_{\operatorname*{reg},A_{\varepsilon}}\Vert f\Vert_{t,\Omega}+\Vert\nabla
u_{\varepsilon,h}\Vert_{t,\Omega}\right)  ^{1-\theta(p,t)}.
\end{align}
Notice that for $f\in L^{t}(\Omega)$ we have $\Vert F \Vert_{-1,t,\Omega}%
\leq\Vert f\Vert_{t,\Omega}$. The combination of \eqref{eq:1} and \eqref{eq:2}
yields the desired estimate. \endproof
Theorem \ref{theo_error} requires several comments. First, from Theorem
\ref{theo_error} it follows that
\begin{multline}
\label{remark1}\Vert\nabla(u-u_{\varepsilon,h})\Vert_{A_{\scriptscriptstyle 0}
}\leq\left\vert \left\vert \left\vert D_{\varepsilon}\right\vert \right\vert
\right\vert _{\infty,\Omega}^{1/2}\mathcal{M}_{\Omega} (u_{\varepsilon
,h},y,\gamma)\\
+\left\vert \left\vert \left\vert B_{\varepsilon}\right\vert \right\vert
\right\vert _{p^{\prime\prime},\Omega}^{1/2}\left(  \mathcal{M}_{\Omega
}(u_{\varepsilon,h},y,\gamma)\right)  ^{\theta(p,t)}\Upsilon(u_{\varepsilon
,h},f,\theta(p,t)) +\left\vert \left\vert \left\vert B_{\varepsilon}
\right\vert \right\vert \right\vert _{p^{\prime\prime},\Omega}^{1/2}%
\Vert\nabla u_{\varepsilon,h}\Vert_{p,\Omega}.
\end{multline}
Here the quantity $\Upsilon(u_{\varepsilon,h},f,\theta(p,t))$ is fully
computable provided that $C_{\operatorname*{reg},A_{\varepsilon}}$ or a
certain upper bound of it is known and norms of $B_{\varepsilon}$ and
$D_{\varepsilon}$ can be computed a priori. Since these matrices have low
order (typically $d=2$ or $3$), the computation of norms is reduced to
well-known algebraic procedures.

It is easy to see that if $A_{\scriptscriptstyle0}=A_{\varepsilon}$, then
$D_{\varepsilon}=I$ and $B_{\varepsilon}=0$. In this case, the second and the
third terms in the right-hand side vanish and the total error is completely
determined by the discretization error encompassed in the first term. Another
limit case arises if $u_{\varepsilon,h}$ coincides with the exact solution of
the problem (\ref{simpl_problem}). Then, the first two terms can be made
arbitrary small and the overall error is determined by the modelling error
encompassed in the last term. In other words, the first two terms can be made
(at least theoretically) arbitrary small if $h$ tends to zero.

Next, it is worth noticing that a somewhat different \textit{modus operandi}
leads to a simpler upper bound of the error. Indeed, in view of (\ref{eq:1})
and Theorem \ref{TheoReg} we have
\[
\Vert\nabla(u-u_{\varepsilon})\Vert_{A_{\scriptscriptstyle0}}\leq\left\vert
\left\vert \left\vert B_{\varepsilon}\right\vert \right\vert \right\vert
_{p^{\prime\prime},\Omega}^{1/2}\Vert\nabla u_{\varepsilon}\Vert_{p,\Omega
}\leq\left\vert \left\vert \left\vert B_{\varepsilon}\right\vert \right\vert
\right\vert _{p^{\prime\prime},\Omega}^{1/2}C_{reg,A_{\varepsilon}}\Vert
f\Vert_{p,\Omega}.
\]
Hence, we obtain another estimate
\begin{equation}
\Vert\nabla(u-u_{\varepsilon,h})\Vert_{A_{\scriptscriptstyle0}}\leq\left\vert
\left\vert \left\vert D_{\varepsilon}\right\vert \right\vert \right\vert
_{\infty,\Omega}^{1/2}\mathcal{M}_{\Omega}(u_{\varepsilon,h},y,\gamma
)+\left\vert \left\vert \left\vert B_{\varepsilon}\right\vert \right\vert
\right\vert _{p^{\prime\prime},\Omega}^{1/2}C_{reg,A_{\varepsilon}}\Vert
f\Vert_{p,\Omega}, \label{remark2}%
\end{equation}
which contains only two terms associated with the discretization and modelling
errors, respectively. Here, the second term associated with the modelling
error depends only on $\varepsilon$. This majorant may be coarser than
(\ref{remark1}), but it allows us to make an a priori estimation of the
modelling error and decide whether or not the problem (\ref{simpl_problem})
could be used for getting an approximation with the desired accuracy $\delta$.
This fact suggests the strategy of solving a simple (smoothened, averaged)
problem (\ref{simpl_problem}) instead of the complicated problem
(\ref{mainproblem}). For example, if we know that the modelling error is
smaller than $\frac{1}{2}\delta$, then we can concentrate on approximations of
$u_{\varepsilon}$ instead of $u$. It is natural to await that in many cases
$u_{\varepsilon}$ will have better regularity properties than $u$. Then,
$u_{\varepsilon,h}$ will converge to $u_{\varepsilon}$ (as $h\rightarrow0$)
much faster than analogous approximations $u_{h}$ in (\ref{mainproblem}) will
converge to $u$ (it is known that this convergence for problems with
complicated coefficients may be arbitrary slow, see \cite{BabOsb}). Hence, we
obtain an efficient method of getting an approximation with the required
accuracy. Moreover, in this way proper reconstructions of the vector function
$y$ can be performed by various \textquotedblleft gradient
recovery\textquotedblright\ or \textquotedblleft gradient
averaging\textquotedblright\ methods, which has been investigated by many
researchers (see, e.g., \cite{Bank2003_1,Bank2003_2,Krizek1984,Zhang2005}).
Thus, for small $h$ the approximation errors can be controlled by the majorant
$\mathcal{M}_{\Omega}$.

Above observations motivate the idea to replace the $L^{\infty}$-norm of
$B_{\varepsilon}$ by some $L^{p}$-norm. We are interested to make the norm
$\left\vert \left\vert \left\vert B_{\varepsilon}\right\vert \right\vert
\right\vert _{p^{\prime\prime},\Omega}$ small and proportional to a certain
positive power of $\varepsilon$. For example, let $A_{\scriptscriptstyle0}%
=\kappa_{\scriptscriptstyle0}I$ and $A_{\varepsilon}=\kappa_{\varepsilon}I$,
where $\kappa_{\scriptscriptstyle0}$ is a jump function and $\kappa
_{\varepsilon}$ is a continuous piecewise affine function approximating this
jump in the $\varepsilon$--strip. Then elementary computations show that
$\left\vert \left\vert \left\vert B_{\varepsilon}\right\vert \right\vert
\right\vert _{p^{\prime\prime},\Omega}\sim\varepsilon^{1/p^{\prime\prime}}$.
Similar a priori analysis is of course possible for more complicated matrices
$A_{\scriptscriptstyle0}$ and $A_{\varepsilon}$.


\section{$\boldsymbol{W^{1,p}}$-Regularity Results for Second Order Elliptic
Problems with Rough Coefficients}

\label{sec:regularity}
Now our goal is to derive $L^{p}\left(  \Omega\right)  $-regularity estimates
for the gradient of $\nabla u$ for some $p>2$. We start from the Poisson
problem (in this case, $A_{\scriptscriptstyle 0} $ is equal to the identity
matrix $I$). Then, we employ perturbation arguments in order to get the
desired estimates for a uniformly elliptic matrix $A_{\scriptscriptstyle 0}
\in L^{\infty}(\Omega,\mathbb{R}_{\operatorname*{sym}}^{d\times d})$. It is
important, that our estimates depend only on the amplitude of jumps in the
coefficients of $A_{\scriptscriptstyle 0} $.

Consider the following problem: For a given $F\in W^{-1,p}\left(
\Omega\right)  $ find $\psi\in W_{0}^{1,p}\left(  \Omega\right)  $ such that%
\begin{equation}
\int_{\Omega}\left\langle \nabla\psi,\nabla v\right\rangle =F\left(  v\right)
\qquad\forall v\in W_{0}^{1,p^{\prime}}\left(  \Omega\right)  .
\label{Poissonproblem}%
\end{equation}

\begin{theorem}
[\cite{Simader1996}]Let $1<p<\infty$. Then, for every $F\in W^{-1,p}\left(
\Omega\right)  $, the problem (\ref{Poissonproblem}) has a unique solution
$\psi\in W_{0}^{1,p}\left(  \Omega\right)  $ which meets the estimate
\begin{equation}
C_{p}^{-1}\left\Vert \nabla\psi\right\Vert _{p, \Omega}\leq\left\Vert
F\right\Vert _{-1,p, \Omega}\leq\left\Vert \nabla\psi\right\Vert _{p,\Omega}
\label{estPoisson}%
\end{equation}
with the \emph{Laplace }$W^{1,p}$\emph{-regularity constant} $C_{p}$ and%
\[
\left\Vert F\right\Vert _{-1,p,\Omega}:=\sup_{\substack{\phi\in W_{0}%
^{1,p^{\prime}}\left(  \Omega\right)  \\\left\Vert \phi\right\Vert
_{1,p^{\prime},\Omega}\leq1}}\left\vert \int_{\Omega}\left\langle \nabla
\psi,\nabla\phi\right\rangle \right\vert .
\]

\end{theorem}

\begin{remark}
The constant $C_{p}$ is independent of $F$ (and $\psi$) but depends on
$\Omega$, $d$ and $p$. We have $C_{2}=1$ and, for $p>2$, $C_{p}$ is
non-decreasing and continuous in $p$ (cf. \cite{Meyers1963}).
\end{remark}

Let $T:=\left\{  \left(  p,P\right)  :2<P<\infty,\ 2\leq p\leq P\right\}  $
and introduce the function
\[
\eta:T\rightarrow\mathbb{R}\qquad\eta\left(  p,P\right)  :=\frac{\frac{1}%
{2}-\frac{1}{p}}{\frac{1}{2}-\frac{1}{P}}%
\]
as well as the function $p^{\ast}:\left[  0,1\right]  \times\left(
2,\infty\right)  \rightarrow\mathbb{R}$%
\[
p^{\ast}\left(  t,P\right)  :=\underset{2\leq p\leq P}{\operatorname{argmax}%
}\left\{  \left(  C_{P}\right)  ^{-\eta\left(  p,P\right)  }\geq1-t\right\}
.
\]

\begin{figure}[h]
\centering
\subfloat{\includegraphics[scale=0.55]{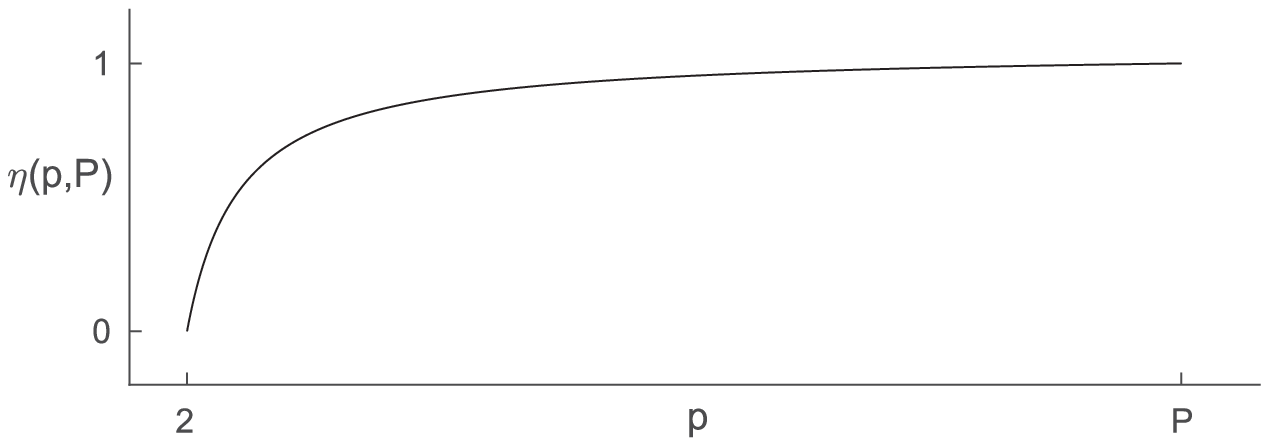}} \hspace{20pt}
\subfloat{\includegraphics[scale=0.55]{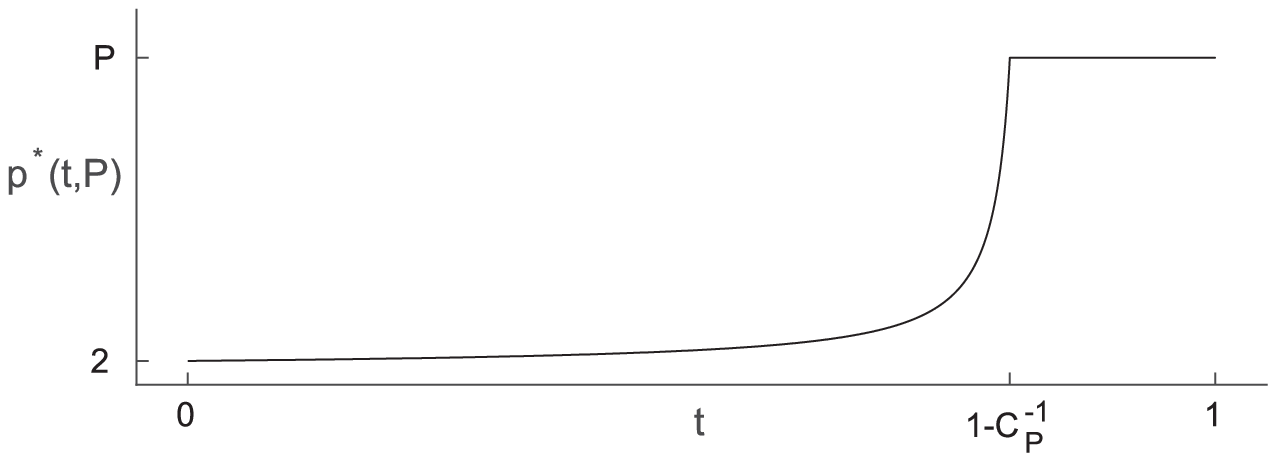}} \caption{The function
$\eta(\cdot,P)$ (left) and the function $p^{*}(\cdot,P)$ (right).}%
\label{fig:p_star}%
\end{figure}

\begin{remark}
For the function $p^{\ast}\left(  \cdot,P\right)  $ it holds%
\begin{equation}
2=p^{\ast}\left(  0,P\right)  \overset{\text{mon. incr.}}{\longrightarrow
}p^{\ast}\left(  1-\frac{1}{C_{P}},P\right)  \overset{\text{for }1-\frac
{1}{C_{P}}\leq t\leq1}{=}p^{\ast}\left(  t,P\right)  =P \label{limitpast}%
\end{equation}
(cf. Figure \ref{fig:p_star}). The function $p^{\ast}$ has the explicit
representation%
\begin{equation}
p^{\ast}\left(  t,P\right)  =\left\{
\begin{array}
[c]{ll}%
\left(  \frac{1}{2}-\dfrac{\log\left(  \dfrac{1}{1-t}\right)  }{\log C_{P}%
}\left(  \dfrac{1}{2}-\dfrac{1}{P}\right)  \right)  ^{-1} & 0\leq t\leq
1-C_{P}^{-1},\\
\quad & \quad\\
P & 1-C_{P}^{-1}<t\leq1.
\end{array}
\right.  \label{defpast2}%
\end{equation}

\end{remark}

Writing $A_{\scriptscriptstyle 0} $ as a perturbation of the identity and
deducing the $L^{p}$-bound of $\nabla u$ from the $L^{p}$-bound for the
solution of the Poisson problem \eqref{Poissonproblem} we obtain the following
result for a uniformly elliptic $L^{\infty}$-coefficient
$A_{\scriptscriptstyle 0} $.

\begin{theorem}
\label{TheoReg}Let $A_{\scriptscriptstyle 0} \in L^{\infty}(\Omega
,\mathbb{R}_{\operatorname*{sym}}^{d\times d})$ satisfy \eqref{defalphabeta},
$F\in W^{-1,P}\left(  \Omega\right)  $ for $2<P<\infty$ and $2\leq p<p^{\ast
}\left(  \frac{\alpha_{\scriptscriptstyle0}}{\beta_{\scriptscriptstyle0}%
},P\right)  $. Then the solution of (\ref{mainproblem}) exists in $W_{0}%
^{1,p}\left(  \Omega\right)  $ and meets the estimate%
\[
\left\Vert \nabla u\right\Vert _{L^{p}\left(  \Omega\right)  }\leq
C_{\operatorname*{reg},A_{\scriptscriptstyle 0} }\left\Vert F\right\Vert
_{W^{-1,p}\left(  \Omega\right)  }%
\]
with
\begin{equation}
C_{\operatorname*{reg},A_{\scriptscriptstyle 0} }:=\frac{1}{\beta
_{\scriptscriptstyle0}}\frac{C_{P}^{\eta\left(  p,P\right)  }}{1-C_{P}%
^{\eta\left(  p,P\right)  }\left(  1-\frac{\alpha_{\scriptscriptstyle0}}%
{\beta_{\scriptscriptstyle0}}\right)  }\text{.} \label{defC}%
\end{equation}

\end{theorem}

For a proof we refer to \cite{Meyers1963, Nochetto2013, Weymuth2016}.

\begin{remark}
\qquad

\begin{enumerate}
\item For the Poisson problem, we have $\alpha_{\scriptscriptstyle0}%
=\beta_{\scriptscriptstyle0}=1$ so that $p^{\ast}\left(  \frac{\alpha
_{\scriptscriptstyle0}}{\beta_{\scriptscriptstyle0}},P\right)  =P$. The
constant in (\ref{defC}) simplifies for $p=P$ to $C_{\operatorname*{reg}%
,A_{\scriptscriptstyle 0} }=C_{P}$ and (\ref{estPoisson}) is reproduced.

\item Let $P\geq2$. For a bounded domain $\Omega$ with $\partial\Omega\in
C^{1}$ there exists some $1\leq C_{\operatorname*{L}}=\mathcal{O}\left(
1\right)  $ such that%
\begin{equation}
C_{P}\leq C_{\operatorname*{L}}P^{d+1} \label{A1}%
\end{equation}
(cf. Remark \ref{rem:constantC_L}).
\end{enumerate}
\end{remark}

\begin{theorem}
Let $\Omega\subset\mathbb{R}^{d}$, $d\geq2$, with $\partial\Omega\in C^{1}$.
Let $\alpha_{\scriptscriptstyle0},\beta_{\scriptscriptstyle0}$ be as in
(\ref{defalphabeta}) and assume that (\ref{A1}) holds.

\begin{enumerate}
\item (small perturbations) For any $P\in\left]  2,\infty\right[  $, consider
$A_{\scriptscriptstyle 0} \in L^{\infty}\left(  \Omega,\mathbb{R}%
_{\operatorname{sym}}^{d\times d}\right)  $ in (\ref{mainproblem}) with
spectral bounds $\alpha_{\scriptscriptstyle0},\beta_{\scriptscriptstyle0}$
such that $\frac{\alpha_{\scriptscriptstyle0}}{\beta_{\scriptscriptstyle0}%
}\geq1-\frac{1}{2C_{\operatorname*{L}}P^{d+1}}$, and let $F\in W^{-1,P}\left(
\Omega\right)  $. Then, for any $2\leq p<P$ it holds%
\[
\left\Vert \nabla u\right\Vert _{L^{p}\left(  \Omega\right)  }\leq
\frac{2C_{\operatorname*{L}}}{\beta_{\scriptscriptstyle0}}P^{d+1}\left\Vert
F\right\Vert _{-1,p,\Omega}.
\]

\item (large perturbations) For any $P\in\left]  2,\infty\right[  $, consider
$A_{\scriptscriptstyle 0} \in L^{\infty}\left(  \Omega,\mathbb{R}%
_{\operatorname{sym}}^{d\times d}\right)  $ in (\ref{mainproblem}) with
spectral bounds $\alpha_{\scriptscriptstyle0},\beta_{\scriptscriptstyle0}$
such that $\frac{\alpha_{\scriptscriptstyle0}}{\beta_{\scriptscriptstyle0}%
}\leq1-C_{P}^{-1}$. For $c\in\left(  0,1\right)  $ define\footnote{By a Taylor
argument it follows that the right-hand side in \eqref{defpmaxest} is bounded
from below by the expression%
\[
p_{\max}\geq2+4c\frac{\frac{\alpha_{\scriptscriptstyle0}}{\beta
_{\scriptscriptstyle0}}}{\log C_{P}}\left(  \frac{1}{2}-\frac{1}{P}\right)
\]
which shows the qualitative dependence on $c$, $\frac{\alpha
_{\scriptscriptstyle0}}{\beta_{\scriptscriptstyle0}}$, $C_{P}$ and $P$
better.}%
\begin{equation}
p_{\max}:=\frac{1}{\frac{1}{2}-c\frac{\log\left(  \frac{1}{1-\frac
{\alpha_{\scriptscriptstyle0}}{\beta_{\scriptscriptstyle0}} }\right)  }{\log
C_{P}}\left(  \frac{1}{2}-\frac{1}{P}\right)  }. \label{defpmaxest}%
\end{equation}
Then, for any $2\leq p<p_{\max}$ it holds%
\[
\left\Vert \nabla u\right\Vert _{L^{p}\left(  \Omega\right)  }\leq\frac
{1}{\beta_{\scriptscriptstyle0}}\frac{1}{\left(  1-\frac{\alpha
_{\scriptscriptstyle0}}{\beta_{\scriptscriptstyle0}}\right)  ^{c}-\left(
1-\frac{\alpha_{\scriptscriptstyle0}}{\beta_{\scriptscriptstyle0}}\right)
}\left\Vert F\right\Vert _{-1,p, \Omega}.
\]

\end{enumerate}
\end{theorem}

\proof

@1 The condition $\frac{\alpha_{\scriptscriptstyle0}}{\beta
_{\scriptscriptstyle0}}\geq1-\frac{1}{2C_{\operatorname*{L}}P^{d+1}} $ implies
$\frac{\alpha_{\scriptscriptstyle0}}{\beta_{\scriptscriptstyle0}}\geq
1-\frac{1}{2C_{P}}> 1-\frac{1}{C_{P}}$ so that (cf. (\ref{limitpast}))%
\[
p^{\ast}\left(  \frac{\alpha_{\scriptscriptstyle0}}{\beta
_{\scriptscriptstyle0}},P\right)  =P.
\]
Hence, we may choose any $2\leq p<P$ for the following. It is easy to see that
then $\eta\left(  p,P\right)  \leq1$ and (\ref{defC}) gives us%
\[
C_{\operatorname*{reg},A_{\scriptscriptstyle 0} }\leq\frac{1}{\beta
_{\scriptscriptstyle0}}\frac{C_{P}}{1-C_{P}\left(  1-\frac{\alpha
_{\scriptscriptstyle0}}{\beta_{\scriptscriptstyle0}}\right)  }\leq\frac
{2}{\beta_{\scriptscriptstyle0}}C_{P}.
\]
This leads to the first assertion.

@2 The condition $0\leq\frac{\alpha_{\scriptscriptstyle0}}{\beta
_{\scriptscriptstyle0}}\leq1-C_{P}^{-1}$ implies that for any $c\in(0,1)$ the
number $p_{\max}$ satisfies%
\[
p_{\max}\leq p^{\ast}\left(  \frac{\alpha_{\scriptscriptstyle0}}%
{\beta_{\scriptscriptstyle0}},P\right)  ,
\]
and Theorem \ref{TheoReg} implies $\nabla u\in L^{p}\left(  \Omega\right)  $
for $2\leq p<p_{\max}$. Note that then $\eta\left(  p,P\right)  \leq
\eta\left(  p_{\max},P\right)  =c\frac{\log\left(  \frac{1}{1-\frac{\alpha
_{0}}{\beta_{\scriptscriptstyle0}} }\right)  }{\log C_{P}}$ so that%
\[
C_{P}^{\eta\left(  p,P\right)  }\leq C_{P}^{\eta\left(  p_{\max},P\right)
}=\left(  \frac{1}{1-\frac{\alpha_{\scriptscriptstyle0}}{\beta
_{\scriptscriptstyle0}}}\right)  ^{c}.
\]
For the regularity constant we obtain%
\[
C_{\operatorname*{reg},A_{\scriptscriptstyle 0} }=\frac{1}{\beta
_{\scriptscriptstyle0}}\frac{1}{\left(  1-\frac{\alpha_{\scriptscriptstyle0}
}{\beta_{\scriptscriptstyle0}}\right)  ^{c}-\left(  1-\frac{\alpha
_{\scriptscriptstyle0}}{\beta_{\scriptscriptstyle0}}\right)  }%
\]
and the assertion follows.\endproof


\section{Analysis of the Laplace $\boldsymbol{W^{1,p}}$-Regularity Constant
for the Full Space Problem}

\label{sec:constant}

Let $f$ be a measurable function on a domain $\Omega$ (bounded or unbounded)
in $\mathbb{R}^{d}$. The distribution function $\mu_{f}$ is defined by
\footnote{For a measurable subset $M\subset\mathbb{R}^{d}$ we set
$|M|:=\int_{M}1$.}
\[
\mu_{f}(t)=\left\vert \left\{  x\in\Omega:f(x)>t\right\}  \right\vert
\]
for $t>0$ and measures the relative size of $f$. A basic property of $\mu$ is
given by the following lemma which is proved in \cite[Lemma 9.7]{Gilbarg1983}.

\begin{lemma}
\label{lemma_distribution_fct} For any $p>0$ and $|f|^{p}\in L^{1}(\Omega)$,
we have
\[
\mu_{f}(t)\leq t^{-p}\int_{\Omega}|f|^{p}.
\]

\end{lemma}

For $f\in L^{p}(\Omega)$, $1<p<\infty$, the Newtonian potential of $f$ is
defined as
\[
\mathcal{N}f(x)=\int_{\Omega}G(x-y)f(y)dy,
\]
where $G$ is the fundamental solution of Laplace's equation which is given by
\[
G(z):=%
\begin{cases}
-\frac{1}{2\pi}\log\Vert z\Vert & d=2,\\
\frac{\Gamma(d/2)}{2\pi^{d/2}(d-2)}\Vert z\Vert^{2-d} & d\geq3,
\end{cases}
\]
where $\Gamma\left(  \cdot\right)  $ denotes the Gamma function. For fixed
$i,j$ we define the linear operator $T:L^{2}(\Omega)\rightarrow L^{2}(\Omega)$
by
\[
Tf:=\partial_{i}\partial_{j}\mathcal{N}f.
\]
Further for $d\geq2$ and $1<p<2$ we define the constants
\begin{equation}
C(d):=2^{d+2}+2^{d+1}d\left(  d+5\right)  +\frac{2\pi^{d/2}}{\Gamma
(d/2)}d^{d/2-1} \label{const_Cd}%
\end{equation}
and
\begin{equation}
C(d,p):=2\left(  \frac{p}{p-1}+\frac{p}{2-p}\right)  ^{1/p}C(d)^{2/p-1}.
\label{const_Cdp}%
\end{equation}

\begin{notation}
For $r>0$ and $x\in\mathbb{R}^{d}$, we denote the ball with radius $r$ around
$x$ by $B_{r}\left(  x\right)  :=\left\{  y\in\mathbb{R}^{d}:\left\Vert
y-x\right\Vert <r\right\}  $. The dimension $d$ is clear from the argument
$x$. We write short $B_{r}$ for $B_{r}\left(  0\right)  $.
\end{notation}

\begin{lemma}
\label{lemma_mu_Tf} Let $\Omega\subset\mathbb{R}^{d}$ be a bounded domain,
$f\in L^{2}(\Omega)$ and $t>0$. Then we have the estimate
\[
\mu_{Tf}(t)\leq C(d)\frac{\Vert f\Vert_{1,\Omega}}{t}%
\]
with $C(d)$ as in \eqref{const_Cd}.
\end{lemma}

\proof

We follow the proof of \cite[Theorem 9.9]{Gilbarg1983} and track the
dependence of the constants on $p$ and $d$. We first extend $f$ to vanish
outside $\Omega$ and fix a cube $K_{0}\supset\Omega$ so that for fixed $t>0$
we have
\[
\int_{K_{0}}f\leq t|K_{0}|.
\]
By bisection of the edges of $K_{0}$, we subdivide $K_{0}$ into $2^{d}$
congruent subcubes with disjoint interiors. Those subcubes $K$ which satisfy
\[
\int_{K}f\leq t|K|
\]
are similarly subdivided and the process is repeated indefinitely. In this way
we obtain a sequence of parallel subcubes $(K_{\ell})_{\ell=1}^{\infty}$ such
that
\begin{equation}
t<\frac{1}{|K_{\ell}|}\int_{K_{\ell}}|f|<2^{d}t \label{eq_subcubes}%
\end{equation}
and
\[
|f|\leq t\quad\text{a.e. on}\ G=K_{0}\backslash\bigcup_{\ell}K_{\ell}.
\]
The function $f$ is now split into a \textquotedblleft good
part\textquotedblright\ $g$ defined by
\[
g(x):=%
\begin{cases}
f(x) & x\in G\\
\frac{1}{|K_{\ell}|}\int_{K_{\ell}}f & x\in K_{\ell}%
\end{cases}
\]
and a \textquotedblleft bad part\textquotedblright\ $b=f-g$. Clearly,
\[
|g|\leq2^{d}t\quad\text{a.e.},\qquad b(x)=0\ \ \text{for}\ x\in G\quad
\text{and}\quad\int_{K_{\ell}}b=0\quad\text{for}\ \ell=1,2,\dots
\]
Since $T$ is linear we have $Tf=Tg+Tb$ and thus
\begin{equation}
\mu_{Tf}(t)\leq\mu_{Tg}(t/2)+\mu_{Tb}(t/2). \label{eq_mu_Tf}%
\end{equation}

From \cite[p. 232]{Gilbarg1983} we know that
\begin{equation}
\mu_{Tg}(t/2)\leq\frac{2^{d+2}}{t}\int|f|. \label{eq_mu_Tg}%
\end{equation}

In a next step we want to estimate $\mu_{Tb}(t/2)$. Writing
\[
b_{\ell}:=b\chi_{K_{\ell}}=%
\begin{cases}
b & \text{on}\ K_{\ell}\\
0 & \text{elsewhere},
\end{cases}
\]
we have
\[
Tb=\sum_{\ell=1}^{\infty}Tb_{\ell}.
\]
Let us now fix some $\ell$ and a sequence $b_{\ell m}\subset C_{0}^{\infty
}(K_{\ell})$ converging to $b_{\ell}$ in $L^{2}(\Omega)$ and satisfying
\[
\int_{K_{\ell}}b_{\ell m}=\int_{K_{\ell}}b_{\ell}=0.
\]
Then for $x\notin K_{\ell}$ we have the relation
\begin{equation}
Tb_{\ell m}(x)=\int_{K_{\ell}}\left(  \partial_{i}\partial_{j}G(x-y)-\partial
_{i}\partial_{j}G(x-\bar{y}_{\ell})\right)  b_{\ell m}(y)dy \label{eq_Tb}%
\end{equation}
where $\bar{y}_{\ell}$ denotes the center of $K_{\ell}$. Further for $x\notin
K_{\ell}$ and $y\in K_{\ell}$ it holds
\[
\partial_{i}\partial_{j}G\left(  x-y\right)  -\partial_{i}\partial_{j}G\left(
x-\bar{y}_{\ell}\right)  =\left.  \left\langle \nabla_{w}\partial_{i}%
\partial_{j}G\left(  x-w\right)  ,\bar{y}_{\ell}-y\right\rangle \right\vert
_{w=\zeta_{y}},
\]
for some point $\zeta_{y}$ between $y$ and $\bar{y}_{\ell}$. Some computations
show that
\[
\left\vert \partial_{i}\partial_{j}\partial_{k}G\left(  z\right)  \right\vert
\leq\frac{d\cdot\Gamma(d/2)}{2\pi^{d/2}}\left\vert \left(  -\frac{\delta
_{ij}+\delta_{ik}+\delta_{jk}}{\left\Vert z\right\Vert ^{d+1}}+\left(
d+2\right)  \frac{1}{\left\Vert z\right\Vert ^{d+1}}\right)  \right\vert ,
\]
where $\delta_{ij}$ denotes the Kronecker symbol. This leads to the estimate
\[
\left\vert \partial_{i}\partial_{j}G\left(  x-y\right)  -\partial_{i}%
\partial_{j}G\left(  x-\bar{y}_{\ell}\right)  \right\vert \leq\frac{d\left(
d+5\right)  \Gamma(d/2)}{2\pi^{d/2}}\cdot\frac{\left\Vert \bar{y}_{\ell
}-y\right\Vert }{\left\Vert x-\zeta_{y}\right\Vert ^{d+1}}.
\]
Note that
\[
\left\Vert x-\zeta_{y}\right\Vert \geq\operatorname*{dist}\left(  x,K_{\ell
}\right)
\]

so that
\begin{equation}
\left\vert \partial_{i}\partial_{j}G\left(  x-y\right)  -\partial_{i}%
\partial_{j}G\left(  x-\bar{y}_{\ell}\right)  \right\vert \leq\frac{d\left(
d+5\right)  \Gamma(d/2)}{2\pi^{d/2}}\cdot\frac{\delta_{\ell}}%
{2\operatorname*{dist}^{d+1}\left(  x,K_{\ell}\right)  }, \label{eq_partialG}%
\end{equation}

where $\delta_{\ell}:=\operatorname*{diam}K_{\ell}$. The combination of
\eqref{eq_Tb} and \eqref{eq_partialG} yields
\[
\left\vert Tb_{\ell m}\left(  x\right)  \right\vert \leq\frac{d\left(
d+5\right)  \Gamma(d/2)\delta_{\ell}}{4\pi^{d/2}\operatorname*{dist}%
^{d+1}(x,K_{\ell})}\int_{K_{\ell}}\left\vert b_{\ell m}\left(  y\right)
\right\vert dy.
\]
It is easy to see that for $x\notin B_{\delta_{\ell}}\left(  \bar{y}_{\ell
}\right)  $ it holds
\[
\operatorname*{dist}\left(  x,K_{\ell}\right)  \geq\operatorname*{dist}\left(
x,B_{\delta_{\ell}/2}\left(  \bar{y}_{\ell}\right)  \right)  \geq\left\Vert
x\right\Vert \inf_{r\geq\delta_{\ell}}\frac{r-\delta_{\ell}/2}{r}=\frac{1}%
{2}\left\Vert x\right\Vert .
\]

Thus,
\begin{align*}
\int_{K_{0}\backslash B_{\delta_{\ell}}\left(  \bar{y}_{\ell}\right)
}\left\vert Tb_{\ell m}\left(  x\right)  \right\vert dx  &  \leq\frac{d\left(
d+5\right)  \Gamma(d/2)\delta_{\ell}}{4\pi^{d/2}}\left(  \int_{\left\Vert
x\right\Vert \geq\delta}\left(  \frac{2}{\left\Vert x\right\Vert }\right)
^{d+1}dx\right)  \int_{K_{\ell}}\left\vert b_{\ell m}\left(  y\right)
\right\vert dy\\
&  =d\left(  d+5\right)  2^{d}\int_{K_{\ell}}\left\vert b_{\ell m}\left(
y\right)  \right\vert dy.
\end{align*}

We set $F^{\ast}:=\bigcup_{\ell}B_{\delta_{\ell}}\left(  \bar{y}_{\ell
}\right)  $ and $G^{\ast}:=K_{0}\backslash F^{\ast}$. Letting $m\rightarrow
\infty$ and summing over $\ell$ we get
\begin{align*}
\int_{G^{\ast}}\left\vert Tb\right\vert  &  \leq\lim_{m\rightarrow\infty}%
\sum_{\ell=1}^{\infty}\int_{K_{0}\backslash B_{\delta_{\ell}}\left(  \bar
{y}_{\ell}\right)  }\left\vert Tb_{\ell m}\right\vert \\
&  \leq d\left(  d+5\right)  2^{d}\lim_{m\rightarrow\infty}\sum_{\ell
=1}^{\infty}\int_{K_{\ell}}\left\vert b_{\ell m}\right\vert \\
&  \leq d\left(  d+5\right)  2^{d}\int\left\vert f\right\vert .
\end{align*}

By Lemma \ref{lemma_distribution_fct} and using the last estimate we obtain%
\begin{equation}
\left\vert \{x\in G^{\ast}:|Tb|>t/2\}\right\vert \leq\left(  \frac{t}%
{2}\right)  ^{-1}\int_{G^{\ast}}\left\vert Tb\right\vert \leq d\left(
d+5\right)  2^{d+1}\frac{\Vert f\Vert_{1,\Omega}}{t}. \label{eq_G_star2}%
\end{equation}

Moreover by \cite[p. 234]{Gilbarg1983} and \eqref{eq_subcubes} we have with
$F:=\bigcup_{\ell}K_{\ell}$
\begin{align}
\label{eq_F_star}\left|  F^{*}\right|  \leq\frac{2\pi^{d/2}}{d\cdot
\Gamma(d/2)}d^{d/2}\left|  F\right|  \leq\frac{2\pi^{d/2}}{d\cdot\Gamma
(d/2)}d^{d/2}\frac{\|f\|_{1,\Omega}}{t}.
\end{align}

Finally we get by \eqref{eq_mu_Tf}, \eqref{eq_mu_Tg}, \eqref{eq_G_star2} and
\eqref{eq_F_star}
\begin{align*}
\mu_{Tf}\left(  t\right)   &  \leq\mu_{Tg}\left(  t/2\right)  +\mu_{Tb}\left(
t/2\right) \\
&  \leq\frac{2^{d+2}}{t}\Vert f\Vert_{1,\Omega}+\left\vert \{x\in G^{\ast
}:|Tb|>t/2\}\right\vert +|F^{\ast}|\\
&  \leq\left(  2^{d+2}+d\left(  d+5\right)  2^{d+1}+\frac{2\pi^{d/2}}%
{d\cdot\Gamma(d/2)}d^{d/2}\right)  \frac{\Vert f\Vert_{1,\Omega}}{t}\\
&  =C(d)\frac{\Vert f\Vert_{1,\Omega}}{t}.
\end{align*}%
\endproof

\endproof

\begin{theorem}
[Calderon-Zygmund estimate]\label{theo:Calderon-Zygmund} Let $\Omega
\subset\mathbb{R}^{d}$ be a bounded domain and $f\in L^{p}(\Omega)$,
$1<p<\infty$. Then it holds
\[
\Vert Tf\Vert_{p,\Omega}\leq C_{1}(d,p)\Vert f\Vert_{p,\Omega}%
\]
with
\begin{equation}
C_{1}(d,p):=%
\begin{cases}
C(d,p) & 1<p\leq\frac{3}{2}\\
C^{\frac{3}{p}(2-p)}\left(  d,\frac{3}{2}\right)  & \frac{3}{2}<p\leq2\\
C^{\frac{3}{p^{\prime}}(2-p^{\prime})}(d,\frac{3}{2}) & 2\leq p< 3\\
C(d,p^{\prime}) & 3\leq p<\infty
\end{cases}
\label{C_1}%
\end{equation}
and $C(d,p)$ as in \eqref{const_Cdp}.
\end{theorem}

\proof

For the case $p=2$ we refer to \cite[Theorem 9.9]{Gilbarg1983}.\newline By
Lemma \ref{lemma_distribution_fct} and since $\Vert Tf\Vert_{2,\Omega}=\Vert
f\Vert_{2,\Omega}$ we know that
\[
\mu_{Tf}(t)\leq\left(  \frac{\Vert f\Vert_{2,\Omega}}{t}\right)  ^{2}%
\]
for all $t>0$ and all $f\in L^{2}(\Omega)$.

Further by Lemma \ref{lemma_mu_Tf} we have
\begin{equation}
\mu_{Tf}(t)\leq C(d)\frac{\Vert f\Vert_{1,\Omega}}{t} \label{eq}%
\end{equation}
for all $t>0$ and all $f\in L^{2}(\Omega)$. Thus it follows by the
Marcinkiewicz interpolation theorem (cf. Theorem \ref{theo_Marcinkiewicz} with
$q=1$ and $r=2$) that
\[
\Vert Tf\Vert_{p,\Omega}\leq C(d,p)\Vert f\Vert_{p,\Omega}\quad\text{for
all}\ 1<p<2
\]
with $C(d,p)$ as in \eqref{const_Cdp}. By a duality argument (cf.
\cite[Theorem 9.9]{Gilbarg1983}) we obtain that
\[
\Vert Tf\Vert_{p,\Omega}\leq C(d,p^{\prime})\Vert f\Vert_{p,\Omega}%
\quad\text{for all}\ 2<p<\infty.
\]

Next we employ the Riesz-Thorin interpolation theorem to remove the singular
behavior of $C\left(  d,p\right)  $ as $p\rightarrow2$. Note that
$T:L^{p}\left(  \Omega\right)  \rightarrow L^{p}\left(  \Omega\right)  $ is
continuous for $p\in\left\{  p_{0},p_{1}\right\}  $ with $p_{0}=3/2$ and
$p_{1}=2$ and $T:L^{p_{0}}\left(  \Omega\right)  +L^{p_{1}}\left(
\Omega\right)  \rightarrow L^{p_{0}}\left(  \Omega\right)  +L^{p_{1}}\left(
\Omega\right)  $ (observe that $L^{p_{0}}\left(  \Omega\right)  +L^{p_{1}%
}\left(  \Omega\right)  =L^{p_{0}}\left(  \Omega\right)  $ since $\Omega$ is
bounded). We know that
\[
\left\Vert Tf\right\Vert _{p,\Omega}\leq C\left(  d,p\right)  \left\Vert
f\right\Vert _{p,\Omega}\quad\text{for}\ p\in\left\{  p_{0},p_{1}\right\}  .
\]
Let%
\begin{equation}
\frac{1}{p}=\frac{1-t}{p_{0}}+\frac{t}{p_{1}}. \label{defpt}%
\end{equation}
Then the Riesz-Thorin interpolation theorem implies that $T:L^{p}\left(
\Omega\right)  \rightarrow L^{p}\left(  \Omega\right)  $ is bounded and
\[
\left\Vert Tf\right\Vert _{p,\Omega}\leq C^{1-t}\left(  d,p_{0}\right)
C^{t}\left(  d,p_{1}\right)  \left\Vert f\right\Vert _{p,\Omega}=C^{\frac
{3}{p}(2-p)}\left(  d,\frac{3}{2}\right)  \Vert f\Vert_{p,\Omega}%
\]
for $\frac{3}{2}<p<2$. Note that $p_{0}^{\prime}=3$ and $p_{1}^{\prime}=2$.
Let
\[
\frac{1}{p}=\frac{1-t}{p_{1}^{\prime}}+\frac{t}{p_{0}^{\prime}}.
\]
Applying again the Riesz-Thorin interpolation theorem yields
\[
\left\Vert Tf\right\Vert _{p,\Omega}\leq C^{1-t}\left(  d,p_{1}^{\prime
}\right)  C^{t}\left(  d,p_{0}^{\prime}\right)  \left\Vert f\right\Vert
_{p,\Omega}=C^{\frac{3}{p^{\prime}}(2-p^{\prime})}\left(  d,\frac{3}%
{2}\right)  \Vert f\Vert_{p,\Omega}%
\]
for $2<p<3$. In total we have proved that
\[
\Vert Tf\Vert_{p,\Omega}\leq C_{1}(d,p)\Vert f\Vert_{p,\Omega}%
\]
with $C_{1}(d,p)$ as in \eqref{C_1}.
\endproof

\endproof

\begin{figure}[h]
\centering
\subfloat{\includegraphics[scale=0.55]{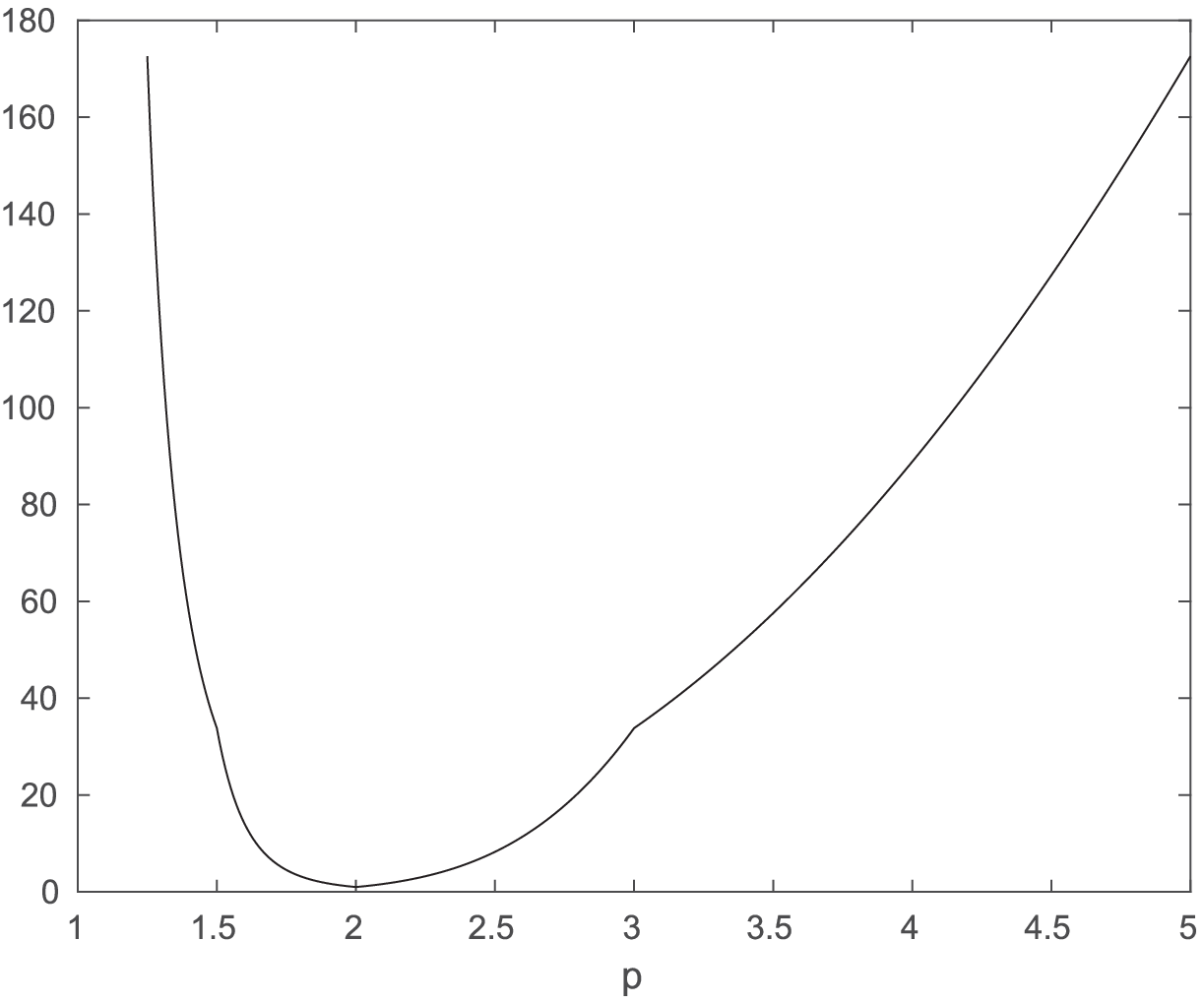}} \hspace{20pt}
\subfloat{\includegraphics[scale=0.55]{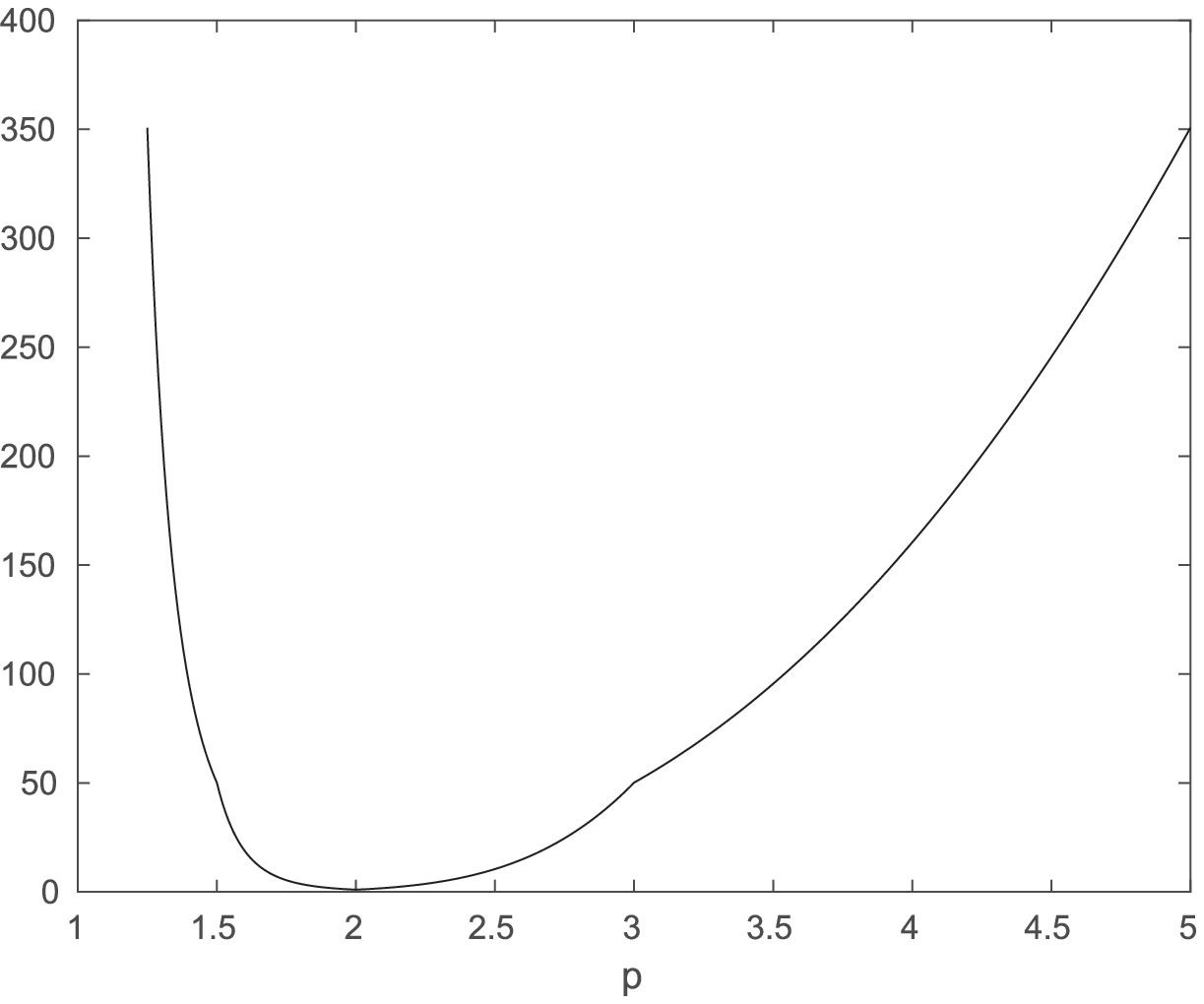}}\caption{The constant
$C_{1}(d,p)$ as a function of $p$ for $d=2$ (left) and for $d=3$ (right).}%
\end{figure}

\begin{remark}
\label{rem_orderC_1} \quad

\begin{enumerate}
\item Note that $C_{1}(d,p)=C_{1}(d,p^{\prime})$ for $1<p<\infty$. Further
observe that $C_{1}(d,p)\geq1$.

\item Since for $p\geq3$ we have
\[
C_{1}(d,p)=2\left(  \frac{p(p-1)}{p-2}\right)  ^{(p-1)/p}C(d)^{1-2/p},
\]
the constant $C_{1}(d,p)$ is of order $\mathcal{O}(p)$ as $p\rightarrow\infty$.
\end{enumerate}
\end{remark}

\begin{remark}
\label{rem:constantC_L} Observe that $T$ can be defined as a bounded operator
on $L^{p}(\Omega)$ even if $\Omega$ is unbounded. In this case Theorem
\ref{theo:Calderon-Zygmund} still holds provided that $d\geq3$ (cf.
\cite{Gilbarg1983}). The constant $C_{1}(d,p)$ in the Calderon-Zygmund
estimate is equal to the Laplace $W^{1,p}$-regularity constant for the full
space problem (cf. \cite[Chapter II, Lem. 2.1]{Simader1996}).\newline The
analysis of the Laplace $W^{1,p}$-regularity constant for bounded domains with
$C^{1}$-boundary can be found in \cite[Chapter II]{Simader1996}. However, the
usual localization techniques such as using cutoff functions on overlapping
balls does not lead to sharp estimates.\newline From the analysis in
\cite{Simader1996} it follows that for a bounded domain $\Omega$ with
$\partial\Omega\in C^{1}$ and $2\leq p<\infty$ there exists some $1\leq
C_{\operatorname*{L}}=\mathcal{O}\left(  1\right)  $ such that%
\[
C_{p}\leq C_{\operatorname*{L}}p^{d+1}%
\]
with $C_{p}$ as in \eqref{estPoisson}.
\end{remark}

\appendix

\section{Interpolation Estimates}

\begin{lemma}
\label{interpolation_estimate} Let $2<r<t<\infty$ and $\theta\in(0,1)$ be such
that
\[
\frac{1}{r}=\frac{\theta}{2}+\frac{1-\theta}{t}.
\]
Then if $u\in L^{t}(\Omega)$, we have the estimate
\begin{equation}
\label{interpolation}\|u\|_{r,\Omega}\leq\|u\|_{2,\Omega}^{\theta
}\|u\|_{t,\Omega}^{1-\theta}.
\end{equation}

\end{lemma}

A proof can be found in \cite{Adams2003}. Note that inequality
\eqref{interpolation} also holds for vector-valued functions.

\begin{theorem}
[Marcinkiewicz interpolation theorem]\label{theo_Marcinkiewicz} Let $S$ be a
linear mapping from $L^{q}(\Omega)\cap L^{r}(\Omega)$ into itself, $1\leq
q<r<\infty$ and suppose that there are constants $S_{1}$ and $S_{2}$ such
that
\[
\mu_{Sf}(t)\leq\left(  \frac{S_{1}\|f\|_{q,\Omega}}{t}\right)  ^{q},\qquad
\mu_{Sf}(t)\leq\left(  \frac{S_{2}\|f\|_{r,\Omega}}{t}\right)  ^{r}%
\]
for all $f\in L^{q}(\Omega)\cap L^{r}(\Omega)$ and $t>0$. Then $S$ extends as
a bounded linear mapping from $L^{p}(\Omega)$ into itself for any $p$ such
that $q<p<r$ and
\[
\|Sf\|_{p,\Omega}\leq2\left(  \frac{p}{p-q}+\frac{p}{r-p}\right)  ^{1/p}
S_{1}^{\alpha}S_{2}^{1-\alpha}\|f\|_{p,\Omega}%
\]
for all $f\in L^{q}(\Omega)\cap L^{r}(\Omega)$, where
\[
\frac{1}{p}=\frac{\alpha}{q}+\frac{1-\alpha}{r}.
\]

\end{theorem}

For a proof we refer to \cite[Theorem 9.8]{Gilbarg1983}.


\begin{thebibliography}{99}                                                                                               %


\bibitem {Adams2003}R. A. Adams and J. J. F. Fournier. Sobolev Spaces, volume
140 of Pure and Applied Mathematics. Elsevier, Amsterdam, 2nd edition, 2003.

\bibitem {Babuska2004}I. Babu\v{s}ka, U. Banerjee, and J. E. Osborn.
Generalized finite element methods -- main ideas, results and perspective.
Int. J. Comput. Methods, 1(1):67--103, 2004.

\bibitem {Babuska1994}I. Babu\v{s}ka, G. Caloz, and J. E. Osborn. Special
finite element methods for a class of second order elliptic problems with
rough coefficients. SIAM J. Numer. Anal., 31(4):945--981, 1994.

\bibitem {Babuska1997}I. Babu\v{s}ka and J. M. Melenk. The partition of unity
method. Int. J. Numer. Meths. Engng., 40(4):727--758, 1997.

\bibitem {BabOsb}I. Babu\v{s}ka and J. E. Osborn. Can a finite element method
perform arbitrarily badly? Math. Comp. 69(230):443--462, 2000.

\bibitem {Bank2003_1}R. E. Bank and J. Xu. Asymptotically exact a posteriori
error estimators, part I: grids with superconvergence. SIAM J. Numer. Anal.,
41(6):2294--2312, 2003.

\bibitem {Bank2003_2}R. E. Bank and J. Xu. Asymptotically exact a posteriori
error estimators, part II: general unstructured grids. SIAM J. Numer. Anal.,
41(6):2313--2332, 2003.

\bibitem {Bensoussan}A. Bensoussan, J.-L. Lions, and G. Papanicolaou.
Asymptotic Analysis for Periodic Structures. North-Holland, Amsterdam, 1978.

\bibitem {Nochetto2013}A. Bonito, R. A. Devore, and R. H. Nochetto. Adaptive
finite element methods for elliptic problems with discontinuous coefficients.
SIAM J. Numer. Anal., 51(6):3106--3134, 2013.

\bibitem {Cioranescu}Doina Cioranescu and Patrizia Donato. An introduction to
homogenization. The Clarendon Press Oxford University Press, New York, 1999.

\bibitem {E_Enquist}Weinan E, Bjorn Engquist, Xiantao Li, Weiqing Ren, and
Eric Vanden-Eijnden. Hetero- geneous multiscale methods: a review. Commun.
Comput. Phys., 2(3):367--450, 2007.

\bibitem {Gilbarg1983}D. Gilbarg and N. S. Trudinger. Elliptic Partial
Differential Equations of Second Order. Springer, Berlin, Heidelberg, New
York, 2nd edition, 1983.

\bibitem {jikov94}V. V. Jikov, S. M. Kozlov, and O. A. Oleinik. Homogenization
of Differential Operators and Integral Functionals. Springer, Berlin, 1994.

\bibitem {Krizek1984}M. Kri\v{z}ek and P. Naittaanm\"aki. Superconvergence
phenomenon in the finite element method arising from averaging of gradients.
Numer. Math., 45:105--116, 1984.

\bibitem {Melenk1996}J. M. Melenk and I. Babu\v{s}ka. The partition of unity
finite element method: Basic theory and applications. Comput. Methods Appl.
Mech. Engrg., 139(1-4):289--314, 1996.

\bibitem {Meyers1963}N. G. Meyers. An $L^{p}$-estimate for the gradient of
solutions of second order elliptic divergence equations. Ann. Scuola Norm.
Sup. Pisa Cl. Sci. (3), 17(3):189--206, 1963.

\bibitem {PreRumSau2011}T. Preusser, M. Rumpf, S. Sauter, and L. O. Schwen. 3D
Composite Finite Elements for Elliptic Boundary Value Problems with
Discontinuous Coefficients. SIAM Journal on Scientific Computing,
33(5):2115--2143, 2011.

\bibitem {Repin2008}S. I. Repin. A Posteriori Error Estimates for Partial
Differential Equations. Walter de Gruyter, Berlin, 2008.

\bibitem {Repin2012}S. I. Repin, T. S. Samrowski, and S. A. Sauter. Combined A
Posteriori Modeling- Discretization Error Estimate for Elliptic Problems with
Complicated Interfaces. ESAIM: Math. Model. Numer. Anal., 46:1389--1405, 2012.

\bibitem {Simader1996}C. G. Simader and H. Sohr. The Dirichlet Problem for the
Laplacian in Bounded and Unbounded Domains: A New Approach to Weak, Strong and
(2+k)-Solutions in Sobolev-Type Spaces, volume 360 of Pitman Research Notes in
Mathematics Series. Addison Wesley Longman, Harlow, Essex, 1996.

\bibitem {Strouboulis2000_1}T. Strouboulis, I. Babu\v{s}ka, and K. Copps. The
design and analysis of the generalized finite element method. Comput. Methods
Appl. Mech. Engrg., 181(1-3):43--69, 2000.

\bibitem {Strouboulis2000_2}T. Strouboulis, K. Copps, and I. Babu\v{s}ka. The
generalized finite element method: an example of its implementation and
illustration of its performance. Int. J. Numer. Meths. Engng.,
47(8):1401--1417, 2000.

\bibitem {Strouboulis2001}T. Strouboulis, K. Copps, and I. Babu\v{s}ka. The
generalized finite element method. Comput. Methods Appl. Mech. Engrg.,
190(32-33):4081--4193, 2001.

\bibitem {Weymuth2016}M. Weymuth. Adaptive Local Basis for Elliptic Problems
with $L^{\infty}$-Coefficients. PhD thesis, University of Zurich, 2016.

\bibitem {Zhang2005}Z. Zhang and A. Naga. A new finite element gradient
recovery method: Superconvergence property. SIAM J. Sci. Comput.,
26:1192--1213, 2005.
\end{thebibliography}
\end{document}